## INTRODUCTORY TOPICS IN BINARY SET FUNCTIONS

Serban E. Vlad
str. Zimbrului, nr.3, bl.PB68, et.2, ap.11, 3700, Oradea, Romania
serbanvlad@excite.com

***Abstract*** *Let* $X \neq \emptyset$ *an arbitrary set and* $\boldsymbol{U} \subset 2^X$ *a non-empty set of subsets. The function* $\mu : \boldsymbol{U} \to \{0,1\}$ *is called binary set function. If* $\mu$ *is countably additive, then it is called a measure. The paper gives some definitions and properties of these functions, its purpose being that of suggesting the reconstruction of the measure theory within this frame, by analogy with* [1], [2].

### 1. Set Rings and Function Rings

1.1 We note with $\boldsymbol{B}_2$ the set $\{0,1\}$, called the *binary Boole* (or *Boolean*) algebra, together with the discrete topology, the order $0 \leq 1$ and the laws: the logical complement '$\overline{\phantom{x}}$', the reunion '$\cup$', the product '$\cdot$', the modulo 2 sum '$\oplus$', the coincidence '$\otimes$':

| $\overline{\phantom{x}}$ | 0 | 1 |
|---|---|---|
| | 1 | 0 |

a)

| $\cup$ | 0 | 1 |
|---|---|---|
| 0 | 0 | 1 |
| 1 | 1 | 1 |

b)

| $\cdot$ | 0 | 1 |
|---|---|---|
| 0 | 0 | 0 |
| 1 | 0 | 1 |

c)

| $\oplus$ | 0 | 1 |
|---|---|---|
| 0 | 0 | 1 |
| 1 | 1 | 0 |

d)

| $\otimes$ | 0 | 1 |
|---|---|---|
| 0 | 1 | 0 |
| 1 | 0 | 1 |

e)

table (1)

1.2 Let $X \neq \emptyset$ be an arbitrary set, that we shall call the *total set*. In the set $2^X$ of the subsets of $X$, the order is given by the inclusion and the laws are: the complementary relative to $X$: '$\overline{\phantom{x}}$', the reunion '$\vee$', the difference '$-$', the intersection '$\wedge$', the symmetrical difference '$\Delta$' and the coincidence '$\Theta$' that is defined like this:

$$A \,\Theta\, B = \overline{A \,\Delta\, B} \tag{1}$$

1.3 **Theorem** Let $\boldsymbol{U} \subset 2^X$ a set of subsets of $X$. The next statements are equivalent:

a) $A, B \in \boldsymbol{U} \Rightarrow A \vee B,\ A - B \in \boldsymbol{U}$

b) $A, B \in \boldsymbol{U} \Rightarrow A \,\Delta\, B,\ A \wedge B \in \boldsymbol{U}$

and the next statements are equivalent too:

c) $A, B \in \boldsymbol{U} \Rightarrow A \wedge B,\ \overline{B - A} \in \boldsymbol{U}$

d) $A, B \in \boldsymbol{U} \Rightarrow A \,\Theta\, B,\ A \vee B \in \boldsymbol{U}$

1.4 **Remark** In the previous theorem, the conditions a), c); b), d) are dual.

1.5 a) The set $\boldsymbol{U}$ that fulfills one of 1.3 a), b) is called *set ring*, or *ring of subsets of* $X$ (*on* $X$ ). N. Bourbaki calls such a set *clan*.

b) Similarly, if $\boldsymbol{U}$ fulfills one of 1.3 c), d), it is called *set ring*, or *ring of subsets of* $X$ (*on* $X$ ), the dual structure of the structure from a).

1.6 **Remark** a) $(\boldsymbol{U}, \Delta, \wedge)$ is really a non-unitary, commutative ring. Its neuter element is $\emptyset$.

b) $(\boldsymbol{U}, \Theta, \vee)$ is itself a non-unitary, commutative ring. Its neuter element is $X$ .

1.7 a) If $X$ belongs to the ring $(\boldsymbol{U}, \Delta, \wedge)$, then $(\boldsymbol{U}, \Delta, \wedge)$ is called a *set algebra*.

b) If $\emptyset$ belongs to the ring $(\boldsymbol{U}, \Theta, \vee)$, then $(\boldsymbol{U}, \Theta, \vee)$ is called a *set algebra* too.

1.8 **Remark** a) The condition that $(\boldsymbol{U}, \Delta, \wedge)$ is a set algebra ( $(\boldsymbol{U}, \Theta, \vee)$ is a set algebra) implies the one that $\boldsymbol{U}$ is a unitary set ring, because if $X \in \boldsymbol{U}$ (if $\emptyset \in \boldsymbol{U}$ ), then it is the unit of the ring.

b) Generally speaking, the unit, if it exists, is given by $\bigvee_{A\in\boldsymbol{U}} A$ (by $\bigwedge_{A\in\boldsymbol{U}} A$ ).

1.9 **Remark** The set algebras are not what is usually meant by the $F$ -algebra structures, where F is a field.

1.10 Let $f : X \to \boldsymbol{B}_2$ a function. Its *support* is by definition the set:

$$supp\ f = \{x \mid x \in X, f(x) = 1\} \quad (1)$$

1.11 If

$$supp\ f = A \quad (1)$$

$f$ will be noted sometimes with $\chi_A$ . This function is called the *characteristic function* of the set $A \subset X$ .

1.12 Let us define for the set ring $(\boldsymbol{U}, \Delta, \wedge)$, respectively for the set ring $(\boldsymbol{U}, \Theta, \vee)$, the set

$$\boldsymbol{U}' = \{f \mid f : X \to \boldsymbol{B}_2, supp\ f \in \boldsymbol{U}\} \quad (1)$$

1.13 $(\boldsymbol{U}', \oplus, \cdot, \cdot)$ and $(\boldsymbol{U}', \otimes, \cup, \cup)$ are $\boldsymbol{B}_2$-algebras, where $'\cdot'$ is the symbol of two laws: the product of the functions and the product of the functions with scalars (both induced from $\boldsymbol{B}_2$), while $'\cup'$ is the dual of $'\cdot'$.

1.14 The associations

$$\boldsymbol{U} \ni A \leftrightarrow \chi_A \in \boldsymbol{U}'$$

are ring isomorphisms. They allow us many times to identify the set rings $\boldsymbol{U} \subset 2^X$ and the function rings $\boldsymbol{U}' \subset \boldsymbol{B}_2^X$ .

## 2. Additive and Countably Additive Set Functions

2.1 **Theorem** Let $\boldsymbol{U} \subset 2^X$ a non-empty family of subsets of $X$ and $\mu : \boldsymbol{U} \to \boldsymbol{B}_2$ a function.

a) If $(\boldsymbol{U}, \Delta, \wedge)$ is a set ring, then the next statements are equivalent:

a.1) $\forall A, B \in \boldsymbol{U}, A \wedge B = \emptyset \Rightarrow \mu(A \vee B) = \mu(A) \oplus \mu(B)$ (1)

a.2) $\forall A, B \in \boldsymbol{U}, \mu(A \,\Delta\, B) = \mu(A) \oplus \mu(B)$ (2)

b) If $(\boldsymbol{U}, \Theta, \vee)$ is a set ring, then the next statements are equivalent:

b.1) $\forall A, B \in \boldsymbol{U}, A \vee B = X \Rightarrow \mu(A \wedge B) = \mu(A) \otimes \mu(B)$ (3)

b.2) $\forall A, B \in \boldsymbol{U}, \mu(A \,\Theta\, B) = \mu(A) \otimes \mu(B)$ (4)

2.2 a) Let $(\boldsymbol{U}, \Delta, \wedge)$ be a set ring. A function $\mu : \boldsymbol{U} \to \boldsymbol{B}_2$ that fulfills one of the equivalent conditions 2.1 a.1), a.2) is called *additive*, or *finitely additive*.

b) In a dual manner, let $(\boldsymbol{U}, \Theta, \vee)$ be a set ring. A function $\mu : \boldsymbol{U} \to \boldsymbol{B}_2$ that fulfills one of the equivalent conditions 2.1 b.1), b.2) is called *additive**, or *finitely additive**.

2.3 The sets of functions $\boldsymbol{U} \to \boldsymbol{B}_2$ which are additive, respectively additive* are noted with $Ad(\boldsymbol{U})$, respectively $Ad^*(\boldsymbol{U})$. They are naturally organized as $\boldsymbol{B}_2$-linear spaces.

2.4 **Theorem** a) Let $\mu \in Ad(\boldsymbol{U})$. For $A, B \in \boldsymbol{U}$, we have:

a.1) $$\mu(\emptyset) = 0 \quad (1)$$

a.2) $$\mu(A - B) = \mu(A) \oplus \mu(A \wedge B) \quad (2)$$

a.3) $$\mu(A \vee B) \oplus \mu(A \wedge B) \oplus \mu(A \,\Delta\, B) = 0 \quad (3)$$

b) If $\mu \in Ad^*(\boldsymbol{U})$, then the next properties are true:

b.1) $$\mu(X) = 1 \quad (4)$$

b.2) $$\mu(\overline{B - A}) = \mu(A) \otimes \mu(A \vee B) \quad (5)$$

b.3) $$\mu(A \wedge B) \otimes \mu(A \vee B) \otimes \mu(A \,\Theta\, B) = 1 \quad (6)$$

where $A, B \in \boldsymbol{U}$.

2.5 Let $a : \boldsymbol{N} \to \boldsymbol{B}_2$,

$$a_n \overset{def}{=} a(n), n \in \boldsymbol{N} \quad (1)$$

a binary sequence. If the support of $a : \{n \mid n \in \boldsymbol{N}, a_n = 1\}$ is a finite set, then the summation modulo 2 has sense:

$$\underset{n \in \boldsymbol{N}}{\Xi} a_n = \begin{cases} 1, \mid supp\ a \mid is\ odd \\ 0, \mid supp\ a \mid is\ even \end{cases} \quad (2)$$

where we have noted with $|\cdot|$ the number of elements of a finite set and where, by definition:

$$|\emptyset| = 0 \quad (3)$$

is even. If the support of $a$ is not finite, then the symbol $\underset{n \in \boldsymbol{N}}{\Xi} a_n$ refers to a divergent series.

2.6 Let $A : \boldsymbol{N} \to \boldsymbol{B}_2$,

$$A_n \overset{def}{=} A(n), n \in \boldsymbol{N} \quad (1)$$

a sequence of sets. If for any $x \in X$ the set $\{n \mid n \in \mathbf{N}, x \in A_n\}$ is finite, then the symmetrical difference has sense:

$$\underset{n \in \mathbf{N}}{\Delta} A_n = \{x \mid x \in X, \mid \{n \mid n \in \mathbf{N}, x \in A_n\} \mid \textit{is odd}\} \tag{2}$$

and if not, the symbol $\underset{n \in \mathbf{N}}{\Delta} A_n$ refers to a divergent series of sets.

2.7 **Theorem** Let $(\mathbf{U}, \Delta, \wedge) \subset 2^X$ be a set ring and $\mu : \mathbf{U} \to \mathbf{B}_2$ a function. The following statements are equivalent:

a) For any sequence of sets $A_n \in \mathbf{U}, n \in \mathbf{N}$, the conditions

a.1) $$n \neq m \Rightarrow A_n \wedge A_m = \emptyset$$

and

a.2) $$\underset{n \in \mathbf{N}}{\vee} A_n \in \mathbf{U}$$

imply

a.3) $\{n \mid n \in \mathbf{N}, \mu(A_n) = 1\}$ is finite

and

a.4) $$\mu(\underset{n \in \mathbf{N}}{\vee} A_n) = \underset{n \in \mathbf{N}}{\Xi} \mu(A_n) \tag{1}$$

b) For any sequence of sets $A_n \in \mathbf{U}, n \in \mathbf{N}$, the conditions

b.1) $\forall x \in X, \{n \mid n \in \mathbf{N}, x \in A_n\}$ is finite

and

b.2) $$\underset{n \in \mathbf{N}}{\Delta} A_n \in \mathbf{U}$$

imply

b.3) $\{n \mid n \in \mathbf{N}, \mu(A_n) = 1\}$ is finite

and

b.4) $$\mu(\underset{n \in \mathbf{N}}{\Delta} A_n) = \underset{n \in \mathbf{N}}{\Xi} \mu(A_n) \tag{2}$$

**Proof** a) $\Rightarrow$ b) Let $A_n \in \mathbf{U}, n \in \mathbf{N}$ so that b.1), b.2) are true under the form:

$$\forall x \in X, \{n \mid n \in \mathbf{N}, x \in A_n\} \in \{0,1\}$$

$$\underset{n \in \mathbf{N}}{\Delta} A_n = \underset{n \in \mathbf{N}}{\vee} A_n \in \mathbf{U} \tag{3}$$

a.1), a.2) being fulfilled, a.3), a.4) are also fulfilled, thus b.3), b.4) are fulfilled.

b) $\Rightarrow$ a) If $A_n \in \mathbf{U}, n \in \mathbf{N}$ satisfies a.1), a.2), then b.1), b.2) are true, thus b.3), b.4) are true resulting that a.3), a.4) are fulfilled.

2.8 a) A function $\mu : \mathbf{U} \to \mathbf{B}_2$ that satisfies one of the equivalent conditions 2.7 a), b) is called *countably additive*, or *measure*.

b) We take in consideration the duals of 2.5, 2.6, 2.7. A function $\mu : \mathbf{U} \to \mathbf{B}_2$ that fulfills one of the duals of the previous equivalent conditions is called *countably additive** or *measure**.

2.9 The sets of countably additive, respectively countably additive* $\mathbf{U} \to \mathbf{B}_2$ functions are noted with $Ad_c(\mathbf{U})$, respectively with $Ad_c^*(\mathbf{U})$.

These sets are $\mathbf{B}_2$-linear spaces.

2.10 The inclusions $Ad_c(\boldsymbol{U}) \subset Ad(\boldsymbol{U})$, $Ad_c^*(\boldsymbol{U}) \subset Ad^*(\boldsymbol{U})$ are easily shown.

2.11 The terminology of additive function, countably additive function and measure is the same if the domain of the function is a $\boldsymbol{B}_2$-algebra $\boldsymbol{U}'$ included in $\boldsymbol{B}_2^X$, via the identification from 1.14.

## 3. Examples

3.1 Let $X \neq \emptyset$ and $\boldsymbol{U} \subset 2^X$ a set ring. The null function $0: \boldsymbol{U} \to \boldsymbol{B}_2$ is a measure; it is the null element of the linear space $Ad_c(\boldsymbol{U})$.

3.2 Suppose that $\mu: \boldsymbol{U} \to \boldsymbol{B}_2$ is a measure and $A \in \boldsymbol{U}$. The function $\mu_1: \boldsymbol{U} \to \boldsymbol{B}_2$ that is defined by:

$$\mu_1(B) = \mu(A \wedge B), B \in \boldsymbol{U} \tag{1}$$

is a measure, called the *restriction of* $\mu$ *at* $A$.

**Proof** Let $A_n \in \boldsymbol{U}, n \in \boldsymbol{N}$ be disjoint two by two with $\underset{n \in \boldsymbol{N}}{\vee} A_n \in \boldsymbol{U}$, resulting that the sets $A \wedge A_n \in \boldsymbol{U}, n \in \boldsymbol{N}$ are disjoint two by two with

$$\underset{n \in \boldsymbol{N}}{\vee}(A \wedge A_n) = A \wedge \underset{n \in \boldsymbol{N}}{\vee} A_n \in \boldsymbol{U} \tag{2}$$

Because $\mu$ is a measure, the set $\{n \mid n \in \boldsymbol{N}, \mu(A \wedge A_n) = 1\}$ is finite and it is true:

$$\mu_1(\underset{n \in \boldsymbol{N}}{\vee} A_n) = \mu(A \wedge \underset{n \in \boldsymbol{N}}{\vee} A_n) = \mu(\underset{n \in \boldsymbol{N}}{\vee}(A \wedge A_n)) = \tag{3}$$

$$= \underset{n \in \boldsymbol{N}}{\Xi} \mu(A \wedge A_n) = \underset{n \in \boldsymbol{N}}{\Xi} \mu_1(A_n)$$

3.3 We fix $x_0 \in X$. The function $\chi^{\{x_0\}}: \boldsymbol{U} \to \boldsymbol{B}_2$ defined by:

$$\chi^{\{x_0\}}(A) = \chi_A(x_0), A \in \boldsymbol{U} \tag{1}$$

is a measure. More general, the sum of these functions is a measure too and this means that to each finite set $H \subset X$ it is associated a function $\chi^H: \boldsymbol{U} \to \boldsymbol{B}_2$ defined in the following way:

$$\chi^H(A) = \underset{x \in H}{\Xi} \chi_A(x), A \in \boldsymbol{U} \tag{2}$$

When $H$ is the empty set, we find the example 3.1.

3.4 $(S_2, \oplus, \bullet, \cdot)$ is the $\boldsymbol{B}_2$-algebra of the binary sequences $x_n \in \boldsymbol{B}_2, n \in \boldsymbol{N}$, where the sum of the sequences $'\oplus'$, the product of the sequences $'\bullet'$ and the product of the sequences with scalars $'\cdot'$ is made coordinatewise. We mention here that the families of sequences $(x_n^p)_n \in S_2, p \in \boldsymbol{N}$ that are disjoint two by two are these that satisfy:

$$p \neq p' \Rightarrow \forall n, x_n^p \cdot x_n^{p'} = 0 \tag{1}$$

Let $k \in \boldsymbol{N}$ and we define $\mu_k: S_2 \to \boldsymbol{B}_2$ by:

$$\mu_k((x_n)) = x_k, (x_n) \in S_2 \tag{2}$$

- the projection of the vector $(x_n)$ of $S_2$ on the $k$-th coordinate. More general, if $H \subset \boldsymbol{N}$ is a finite set

$$H = \{k_1,...,k_p\} \tag{3}$$

then we have the sum of functions $\mu_H : S_2 \to \boldsymbol{B}_2$,

$$\mu_H = \mu_{k_1} \oplus ... \oplus \mu_{k_p} \tag{4}$$

$\mu_k$ and $\mu_H$ are countably additive; if $H$ is empty, then $\mu_H$ is by definition the null function.

3.5 a) We say that the sequence $x_n \in \boldsymbol{B}_2, n \in \boldsymbol{N}$ *converges to* $x^0 \in \boldsymbol{B}_2$ if

$$\exists N \in \boldsymbol{N}, \forall n \geq N, x_n = x^0 \tag{1}$$

If so, the unique $x^0$ with this property (because $x$ is a function) is called the *limit of* $(x_n)$. If the previous statement is made under the weaker form: the sequence $(x_n)$ is *convergent*, this means that such an $x^0$ like at (1) (uniquely) exists. The limit of the sequence $(x_n)$ has the usual notation $\lim\limits_{n \to \infty} x_n$.

b) $(S_2^0, \oplus, \bullet, \cdot)$ is the $\boldsymbol{B}_2$-algebra of the binary sequences $x_n \in \boldsymbol{B}_2, n \in \boldsymbol{N}$ that converge to $0$. We define the measure $\mu : S_2^0 \to \boldsymbol{B}_2$ by

$$\mu((x_n)) = \underset{n \in \boldsymbol{N}}{\Xi} x_n, (x_n) \in S_2^0 \tag{1}$$

3.6 $(S_2^c, \oplus, \bullet, \cdot)$ is the $\boldsymbol{B}_2$-algebra of the convergent binary sequences $x_n \in \boldsymbol{B}_2, n \in \boldsymbol{N}$ and we define $\mu : S_2^c \to \boldsymbol{B}_2$ by:

$$\mu((x_n)) = \lim_{n \to \infty} x_n, (x_n) \in S_2^c \tag{1}$$

$\mu$ is additive, but it is not countably additive. In order to see this, we give the example of the sequence of convergent sequences (the canonical base of $S_2^c$):

$$\varepsilon^n : \boldsymbol{N} \to \boldsymbol{B}_2, \varepsilon^n(m) = \begin{cases} 1, n = m \\ 0, else \end{cases}, m, n \in \boldsymbol{N} \tag{2}$$

$(\varepsilon^n)_n$ are disjoint two by two, their reunion is the constant 1 sequence that is convergent and on the other hand

$$\mu(\bigcup_{n \in \boldsymbol{N}} \varepsilon^n) = 1 \neq 0 = \underset{n \in \boldsymbol{N}}{\Xi} \mu(\varepsilon^n) \tag{3}$$

3.7 A variant of 3.4 is obtained if we take instead of $S_2 = \boldsymbol{B}_2^{\boldsymbol{N}}$ an arbitrary function $\boldsymbol{B}_2$-algebra $\boldsymbol{U} \subset \boldsymbol{B}_2^X$. Let $x_0 \in X$; the function $\mu_{x_0} : \boldsymbol{U} \to \boldsymbol{B}_2$ defined like this:

$$\mu_{x_0}(f) = f(x_0), f \in \boldsymbol{U} \tag{1}$$

is a measure. More general, if $H \subset X$ is a finite set, then the function $\mu_H : \boldsymbol{U} \to \boldsymbol{B}_2$ defined in the following manner:

$$\mu_H(f) = \underset{x\in H}{\Xi} f(x),\ f \in \boldsymbol{U} \tag{2}$$

is a measure. If $H$ is empty, then by definition $\mu_H$ is the null function.

We mention the fact that $f^p \in \boldsymbol{U},\ p \in \boldsymbol{N}$ are disjoint two by two if

$$p \neq p' \Rightarrow \forall x \in X,\ f^p(x) \cdot f^{p'}(x) = 0 \tag{3}$$

3.8 We note with $R_f(X)$ the ring - relative to $\Delta, \wedge$ - of the finite subsets of $X$. The function $\mu_f^X : R_f(X) \to \boldsymbol{B}_2$,

$$\mu_f^X(A) = \begin{cases} 1, |A|\ is\ odd \\ 0, |A|\ is\ even \end{cases},\ A \in R_f(X) \tag{1}$$

is a measure, called the *finite Boolean measure*.

3.9 We note with $Inf_f$ the ring of the *inferiorly finite sets* $A \subset \boldsymbol{R}$, i.e. the sets with the following property:

$$\forall \alpha \in \boldsymbol{R}, (-\infty, \alpha) \wedge A \text{ is finite}$$

We fix some $\alpha \in \boldsymbol{R}$ and we define $\mu_\alpha : Inf_f \to \boldsymbol{B}_2$ by:

$$\mu_\alpha(A) = \mu_f((-\infty, \alpha) \wedge A),\ A \in Inf_f \tag{1}$$

$\mu_\alpha$ is countably additive: for any family $A_n \in Inf_f, n \in \boldsymbol{N}$ of two by two disjoint sets so that $\underset{n\in \boldsymbol{N}}{\vee} A_n \in Inf_f$, only a finite number of sets $A_n$ fulfill $(-\infty, \alpha) \wedge A_n \neq \emptyset$ etc.

3.10 a) Let $X \subset \boldsymbol{R}$ and $t \in \boldsymbol{R} \vee \{\infty\}$ a point so that

$$\forall t' < t, (t', t) \wedge X \text{ is infinite}$$

b) We say that the function $f : X \to \boldsymbol{B}_2$ has a *left limit in* $t$, noted with $f(t-0) \in \boldsymbol{B}_2$, if the next property is true:

$$\exists t' < t, \forall \xi \in (t', t) \wedge X,\ f(\xi) = f(t-0) \tag{1}$$

c) We note with $Lim_X^-(t)$ the $\boldsymbol{B}_2$-algebra of the $X \to \boldsymbol{B}_2$ functions that have a left limit in $t$.

d) The function $\mu : Lim_X^-(t) \to \boldsymbol{B}_2$;

$$\mu(f) = f(t-0),\ f \in Lim_X^-(t) \tag{1}$$

is a measure, this example being analogue to 3.7.

e) Other examples of measures of the same type with this one may be given.

3.11 a) For $a, b \in \boldsymbol{R} \vee \{\infty\}$, the *symmetrical interval* $[[a,b))$ is defined by:

$$[[a,b)) = \begin{cases} [a,b),\ a < b \\ [b,a),\ b < a \\ \ \emptyset\ \ ,\ b = a \end{cases} \tag{1}$$

b) We note with $Sym^-$ the set ring - relative to $\Delta, \wedge$ - generated by the symmetrical intervals $[[a,b))$.

c) We define $\mu : Sym^- \to \boldsymbol{B}_2$ by:

$$\mu(A)=\begin{cases}1, \; if \; \sup A=\infty\\ 0, \; else\end{cases} \tag{2}$$

where $A\in Sym^{-}$. Because in a sequence of sets $A_n\in Sym^{-}$, $n\in \boldsymbol{N}$ that are disjoint two by two with $\bigvee_{n\in \boldsymbol{N}} A_n\in Sym^{-}$ at most one satisfies the condition $\sup A_n=\infty$, it may be shown that $\mu$ is a measure.

3.12 a) We define the next $\boldsymbol{B}_2$-algebras of functions $f:\boldsymbol{R}\to \boldsymbol{B}_2$:

$$I_{[[a,b))}=\{f \mid [[a,b)) \wedge supp \; f \; is \; finite\}, a,b\in \boldsymbol{R}\vee\{\infty\} \tag{1}$$

$$I_{\infty}=\{f \mid supp \; f \; is \; finite\} \tag{2}$$

and the integrals

$$\int_a^{b} {}^{-} f = \underset{t\in [[a,b))}{\Xi} f(t),\; f\in I_{[[a,b))} \tag{3}$$

$$\int_{-\infty}^{\infty} f = \underset{t\in \boldsymbol{R}}{\Xi} f(t),\; f\in I_{\infty} \tag{4}$$

b) The next $I_{[[a,b))}\to \boldsymbol{B}_2$, $I_{\infty}\to \boldsymbol{B}_2$ functions:

$$\mu(f)=\int_a^{b} {}^{-} f,\; f\in I_{[[a,b))} \tag{5}$$

$$\mu(f)=\int_{-\infty}^{\infty} f,\; f\in I_{\infty} \tag{6}$$

are measures.

3.13 a) The set $\boldsymbol{S}\subset 2^{\boldsymbol{R}}$ defined in the next way:

$$\boldsymbol{S}=\{(a_1,b_1)\,\Delta\,...\,\Delta\,(a_p,b_p)\,\Delta\,\{c_1,...,c_n\} \mid a_1,b_1, \;...$$

$$...,a_p,b_p,c_1,...,c_n\in \boldsymbol{R},\, p,n\in \boldsymbol{N}\} \tag{1}$$

is a ring of subsets of $\boldsymbol{R}$ and we have supposed that

$$p=0 \Rightarrow (a_1,b_1)\,\Delta\,...\,\Delta\,(a_p,b_p)=\varnothing \tag{2}$$

$$n=0 \Rightarrow \{c_1,...,c_n\}=\varnothing \tag{3}$$

b) The function $\mu:\boldsymbol{S}\to \boldsymbol{B}_2$ given by:

$$\mu\{(a_1,b_1)\,\Delta\,...\,\Delta\,(a_p,b_p)\,\Delta\,\{c_1,...,c_n\}\}=\pi(p+n) \tag{4}$$

where $\pi:\boldsymbol{N}\to \boldsymbol{B}_2$ is the *parity function*:

$$\pi(m)=\begin{cases}1, \; if \; m \; is \; odd\\ 0, \; if \; m \; is \; even\end{cases}, m\in \boldsymbol{N} \tag{5}$$

- is additive, but it is not countably additive. In order to see this fact, we take the sequence

$$[\frac{1}{n+2},\frac{1}{n+1})=(\frac{1}{n+2},\frac{1}{n+1})\,\Delta\,\{\frac{1}{n+2}\}\in \boldsymbol{S}, n\in \boldsymbol{N} \tag{6}$$

of sets that are disjoint two by two, satisfying

$$\bigvee_{n\in \boldsymbol{N}} [\frac{1}{n+2},\frac{1}{n+1})=(0,1)\in \boldsymbol{S} \tag{7}$$

$$\{n \mid \mu([\frac{1}{n+2},\frac{1}{n+1}))=1\}=\emptyset \tag{8}$$

$$\mu((0,1))=1\neq 0=\underset{n\in \boldsymbol{N}}{\Xi}\mu([\frac{1}{n+2},\frac{1}{n+1})) \tag{9}$$

3.14 a) We note with

$$R_f^*(X)=\{H \mid H\subset X, \overline{H} \ is \ finite\} \tag{1}$$

This set is a set ring relative to $\Theta, \vee$ and it is the dual structure of $R_f(X)$.

b) A typical example of measure* is given by the function $\mu_f^{*X}: R_f^*(X)\to \boldsymbol{B}_2$ that is defined in the next manner:

$$\mu_f^{*X}(H)=\begin{cases}0, |\overline{H}| \ is \ odd\\ 1, |\overline{H}| \ is \ even\end{cases}=\overline{\mu_f^X(\overline{H})}, H\in R_f^*(X) \tag{2}$$

(In the equations (1), (2) the superior bar notes two things: the complementary of a set and the logical complement.)

Let the sequence of sets $A_n\in R_f^*(X), n\in \boldsymbol{N}$ that are disjoint* two by two:

$$n,m\in \boldsymbol{N}, n\neq m\Rightarrow A_n\vee A_m=X \text{ (i.e. } \overline{A_n}\wedge\overline{A_m}=\emptyset) \tag{3}$$

so that $\underset{n\in \boldsymbol{N}}{\wedge} A_n\in R_f^*(X)$. Because from the definition of $R_f^*(X)$, the set

$$\overline{\underset{n\in \boldsymbol{N}}{\wedge} A_n}=\underset{n\in \boldsymbol{N}}{\vee}\overline{A_n} \tag{4}$$

is finite, there results the existence of a rank $N$ with the property that $\overline{A_n}$ are empty for $n>N$. We have:

$$\mu_f^{*X}(\underset{n\in \boldsymbol{N}}{\wedge} A_n)=\overline{\mu_f^X(\overline{\underset{n\in \boldsymbol{N}}{\wedge} A_n})}=\overline{\mu_f^X(\underset{n\in \boldsymbol{N}}{\vee}\overline{A_n})}=\overline{\mu_f^X(\overline{A_0}\vee\overline{A_1}\vee...\vee\overline{A_N})}=$$

$$=\overline{\mu_f^X(\overline{A_0})\oplus\mu_f^X(\overline{A_1})\oplus...\oplus\mu_f^X(\overline{A_N})}=\overline{\underset{n\in \boldsymbol{N}}{\Xi}\mu_f^X(\overline{A_n})}=\overline{\underset{n\in \boldsymbol{N}}{\Xi}\overline{\overline{\mu_f^X(\overline{A_n})}}}=$$

$$=\overline{\underset{n\in \boldsymbol{N}}{\Xi}\overline{\mu_f^{*X}(A_n)}}=\underset{n\in \boldsymbol{N}}{\otimes}\mu_f^{*X}(A_n) \tag{5}$$

## 4. The Behavior of the Measures Relative to the Monotonous Sequences of Sets

4.1 a) The family $A_n\subset X, n\in \boldsymbol{N}$ is called *ascending sequence* of sets if

$$A_0\subset A_1\subset A_2\subset ... \tag{1}$$

In this case, the reunion $\underset{n\in \boldsymbol{N}}{\vee} A_n$ is called the *limit* of the sequence and is noted sometimes with $\lim\limits_{n\to\infty} A_n$.

b) The family $A_n\subset X, n\in \boldsymbol{N}$ is called *descending sequence* of sets if

$$A_0\supset A_1\supset A_2\supset ... \tag{2}$$

The intersection $\bigwedge_{n\in \boldsymbol{N}} A_n$ is called the *limit* of the sequence and is noted sometimes with $\lim_{n\to\infty} A_n$.

c) If the sequence $A_n \subset X, n\in \boldsymbol{N}$ is either ascending, or descending, then we say that it is *monotonous*.

4.2 **Theorem** Let $\boldsymbol{U} \subset 2^X$ a set ring and the function $\mu : \boldsymbol{U} \to \boldsymbol{B}_2$.

a) Let $A_n \in \boldsymbol{U}, n\in \boldsymbol{N}$ an arbitrary ascending sequence of sets satisfying the property that the set

$$A = \bigvee_{n\in \boldsymbol{N}} A_n \tag{1}$$

belongs to $\boldsymbol{U}$. If $\mu$ is a measure, then the binary sequence $(\mu(A_n))_n$ is convergent (see 3.5 a)) and it is true:

$$\mu(A) = \lim_{n\to\infty} \mu(A_n) \tag{2}$$

b) Suppose that $\mu$ is additive and it satisfies the property: for any ascending sequence $A_n \in \boldsymbol{U}, n\in \boldsymbol{N}$ of sets so that its reunion $A$ belongs to $\boldsymbol{U}$, the binary sequence $(\mu(A_n))_n$ is convergent and the relation (2) takes place. Then $\mu$ is a measure.

**Proof** a) We have the disjoint reunion:

$$A = A_0 \vee (A_1 - A_0) \vee ... \vee (A_{n+1} - A_n) \vee ... \tag{3}$$

Because $\mu$ is a measure, it results that there exists $N \in \boldsymbol{N}$ so that

$$n > N \Rightarrow \mu(A_{n+1} - A_n) = 0 \tag{4}$$

thus

$$\mu(A) = \mu(A_0) \oplus \mu(A_1 - A_0) \oplus ... \oplus \mu(A_{N+1} - A_N) = \tag{5}$$

$$= \mu(A_0) \oplus \mu(A_1) \oplus \mu(A_0) \oplus ... \oplus \mu(A_{N+1}) \oplus \mu(A_N) = \mu(A_{N+1})$$

(4) is equivalent with the convergence of the sequence $(\mu(A_n))_n$, as it can be rewritten under the form:

$$n > N \Rightarrow \mu(A_{n+1}) = \mu(A_n) \tag{6}$$

and (5) is equivalent in this situation with (2). In the last equations, we have used 2.4 a.2) under the form:

$$\mu(A_{n+1} - A_n) = \mu(A_{n+1}) \oplus \mu(A_{n+1} \wedge A_n) = \mu(A_{n+1}) \oplus \mu(A_n), n\in \boldsymbol{N} \tag{7}$$

b) Let $A_n^{'} \in \boldsymbol{U}, n\in \boldsymbol{N}$ a sequence of sets that are disjoint two by two and let us suppose that their reunion

$$A = \bigvee_{n\in \boldsymbol{N}} A_n^{'} \tag{8}$$

belongs to $\boldsymbol{U}$. We define the sequence $A_n \in \boldsymbol{U}, n\in \boldsymbol{N}$ by:

$$A_n = A_0^{'} \vee A_1^{'} \vee ... \vee A_n^{'}, n\in \boldsymbol{N} \tag{9}$$

and it is remarked that it is ascending and (1) is satisfied. The hypothesis states the convergence of the sequence with the general term

$$\mu(A_n) = \mu(A_0^{'}) \oplus \mu(A_1^{'}) \oplus ... \oplus \mu(A_n^{'}) \tag{10}$$

in other words there exists an $N \in \boldsymbol{N}$ for which the implication

$$n > N \Rightarrow \mu(A_n^{'}) = 0 \tag{11}$$

is true. The relation (2) becomes

$$\mu(A)=\mu(A_0^{'})\oplus\mu(A_1^{'})\oplus...\oplus\mu(A_N^{'})=\underset{n\in \mathbf{N}}{\Xi}\mu(A_n^{'}) \tag{12}$$

i.e. $\mu$ is a measure.

4.3 **Theorem** It is considered the set ring $\boldsymbol{U}$ and the function $\mu:\boldsymbol{U}\to \boldsymbol{B}_2$.

a) We suppose that $A_n\in \boldsymbol{U}, n\in \mathbf{N}$ is an arbitrary descending sequence of sets whose intersection

$$A=\underset{n\in \mathbf{N}}{\wedge} A_n \tag{1}$$

belongs to $\boldsymbol{U}$ and that $\mu$ is a measure. Then the binary sequence $(\mu(A_n))_n$ is convergent and it is true:

$$\mu(A)=\lim_{n\to\infty}\mu(A_n) \tag{2}$$

b) Let us suppose that $\mu$ is additive and the next property is satisfied: for any descending sequence $A_n\in \boldsymbol{U}, n\in \mathbf{N}$ of sets so that its intersection $A$ belongs to $\boldsymbol{U}$, the binary sequence $(\mu(A_n))_n$ is convergent and the relation (2) is true. In these circumstances $\mu$ is a measure.

**Proof** a) Let us remark for the beginning that the set

$$\underset{n\in \mathbf{N}}{\vee}(A_0-A_n)=A_0-\underset{n\in \mathbf{N}}{\wedge}A_n=A_0-A \tag{3}$$

belongs to $\boldsymbol{U}$ and the sequence of sets

$$A_0-A_0\subset A_0-A_1\subset A_0-A_2\subset ... \tag{4}$$

is ascending. We apply 4.2 a) resulting that the binary sequence $(\mu(A_0-A_n))_n$ is convergent and that it takes place

$$\mu(A_0-A)=\lim_{n\to\infty}\mu(A_0-A_n) \tag{5}$$

From (5) it results that

$$\mu(A_0)\oplus\mu(A)=\mu(A_0)\oplus\lim_{n\to\infty}\mu(A_n) \tag{6}$$

and we have the validity of (2).

b) Let $A_n^{'}\in \boldsymbol{U}, n\in \mathbf{N}$ a sequence of sets that are disjoint two by two with the property that their reunion

$$A'=\underset{n\in \mathbf{N}}{\vee}A_n^{'} \tag{7}$$

belongs to $\boldsymbol{U}$. We define the sequence of sets from $\boldsymbol{U}$:

$$A_n=A'-(A_0^{'}\vee A_1^{'}\vee ...\vee A_n^{'})=(A'-A_0^{'})\wedge(A'-A_1^{'})\wedge ...\wedge(A'-A_n^{'}) \tag{8}$$

where $n\in \mathbf{N}$ that proves to be descending and its meet

$$A=\underset{n\in \mathbf{N}}{\wedge}A_n=\underset{n\in \mathbf{N}}{\wedge}\overset{n}{\underset{k=0}{\wedge}}(A'-A_k^{'})=A'-\underset{n\in \mathbf{N}}{\vee}\overset{n}{\underset{k=0}{\vee}}A_k^{'}=A'-A'=\varnothing \tag{9}$$

belongs to $\boldsymbol{U}$. The hypothesis states that the binary sequence $(\mu(A_n))_n$ is convergent and the relation (2) becomes:

$$0=\mu(\underset{n\in \mathbf{N}}{\wedge}A_n)=\lim_{n\to\infty}\mu(A_n) \tag{10}$$

There exists a rank $N\in \mathbf{N}$ so that for any $n>N$ we have:

$$0=\mu(A_n)=\mu(A'-(A_0^{'}\vee A_1^{'}\vee ...\vee A_n^{'}))=\mu(A')\oplus\mu(A_0^{'})\oplus\mu(A_1^{'})\oplus ...\oplus\mu(A_n^{'}) \quad (11)$$

We have that $(\mu(A_n^{'}))_n$ converges to $0$ and if $k>N$ then

$$\mu(\underset{n\in \boldsymbol{N}}{\vee} A_n^{'})=\mu(A')=\overset{k}{\underset{n=0}{\Xi}}\mu(A_n^{'})=\underset{n\in \boldsymbol{N}}{\Xi}\mu(A_n^{'}) \quad (12)$$

## 5. Derivable Measures

5.1 In this paragraph we shall consider that the total space $X$ is equal with $\boldsymbol{R}^n$, $n\geq 1$. The elements $x\in X$ will be consequently $n-$tuples $(x_1,...,x_n)\in \boldsymbol{R}^n$.

5.2 We define the family

$$\boldsymbol{U}_n=\{A\,|\,A\subset \boldsymbol{R}^n,\ A\ is\ bounded\} \quad (1)$$

It is a set ring (relative to $\Delta$ and $\wedge$).

5.3 Let $A\in \boldsymbol{U}_n$ be a bounded set. Its *diameter* is defined to be the real non-negative number

$$d(A)=\sup_{x,y\in A}\sqrt{(x_1-y_1)^2+...+(x_n-y_n)^2} \quad (1)$$

5.4 We define the *locally finite* sets from $\boldsymbol{R}^n$ to be these sets $H\subset \boldsymbol{R}^n$ with the property that

$$\forall A\in \boldsymbol{U}_n,\ A\wedge H\ is\ finite$$

5.5 The set of the locally finite sets from $\boldsymbol{R}^n$ is noted with $Loc_f^{(n)}$ and it is a set ring.

5.6 **Proposition** Let us take a set $H\in Loc_f^{(n)}$. The function $\mu_H:\boldsymbol{U}_n\to \boldsymbol{B}_2$ defined by:

$$\mu_H(A)=\pi(|A\wedge H|),\ A\in \boldsymbol{U}_n \quad (1)$$

is a measure (the function $\pi$ was defined at 3.13 (5)).

**Proof** Let $A_p\in \boldsymbol{U}_n,\ p\in \boldsymbol{N}$ a family of sets that are disjoint two by two with the property that $\underset{p\in \boldsymbol{N}}{\vee} A_p\in \boldsymbol{U}_n$. Because $\underset{p\in \boldsymbol{N}}{\vee} A_p\wedge H$ is a finite set, there exists a number $N\in \boldsymbol{N}$ with:

$$p>N\Rightarrow A_p\wedge H=\emptyset \quad (2)$$

We infer that

$$\{p\,|\,\mu_H(A_p)=1\}\subset\{0,1,...,N\} \quad (3)$$

$$\mu_H(\underset{p\in \boldsymbol{N}}{\vee} A_p)=\pi(|\underset{p\in \boldsymbol{N}}{\vee} A_p\wedge H|)=\pi(|\underset{p\in \boldsymbol{N}}{\vee}(A_p\wedge H)|)= \quad (4)$$

$$=\pi(|(A_0\wedge H)\vee(A_1\wedge H)\vee ...\vee(A_N\wedge H)|)=$$

$$=\pi(|A_0\wedge H|+|A_1\wedge H|+...+|A_N\wedge H|)=$$

$$=\pi(|A_0\wedge H|)\oplus\pi(|A_1\wedge H|)\oplus ...\oplus\pi(|A_N\wedge H|)=$$

$$= \underset{p\in \boldsymbol{N}}{\Xi} \pi(| A_p \wedge H |) = \underset{p\in \boldsymbol{N}}{\Xi} \mu_H(A_p)$$

5.7 **Proposition** The function $\mu_H \in Ad_c(\boldsymbol{U}_n)$ that was previously defined fulfills the property that for any $A\in \boldsymbol{U}_n$ and $x\in A$:

$$\exists\varepsilon > 0, \exists a\in \boldsymbol{B}_2, \forall B\in \boldsymbol{U}_n, (x\in B \text{ and } d(B) < \varepsilon) \Rightarrow \mu_H(B) = a \tag{1}$$

**Proof** We define the real positive number

$$\varepsilon = \begin{cases} \min \{\sqrt{(x_1 - y_1)^2 + ... + (x_n - y_n)^2} \mid y\in H - \{x\}\}, & H - \{x\} \neq \varnothing \\ 1 & , H - \{x\} = \varnothing \end{cases} \tag{2}$$

Such an $\varepsilon$ exists, if not there would exist a sphere $S_x$ with the center in $x$ and the property that $S_x \wedge H$ is infinite and this is a contradiction with the hypothesis $H \in Loc_f^{(n)}$.

Any bounded set $B\in \boldsymbol{U}_n$ with the properties that $x\in B$ and $d(B) < \varepsilon$ fulfills the relations:

$$(B - \{x\}) \wedge H = \varnothing \tag{3}$$

$$\mu_H(B - \{x\}) = 0 \tag{4}$$

$$\mu_H(B) = \mu_H((B - \{x\}) \vee \{x\}) = \mu_H(B - \{x\}) \oplus \mu_H(\{x\}) = \tag{5}$$

$$= \mu_H(\{x\}) = \pi(| \{x\} \wedge H |) = \begin{cases} 1, & x\in H \\ 0, & x\notin H \end{cases}$$

5.8 Let now $\mu : \boldsymbol{U}_n \to \boldsymbol{B}_2$ be a measure.

a) We say that it is *derivable in* $x\in A$, where $A\in \boldsymbol{U}_n$, if

$$\exists\varepsilon > 0, \exists a\in \boldsymbol{B}_2, \forall B\in \boldsymbol{U}_n, (x\in B \text{ and } d(B) < \varepsilon) \Rightarrow \mu(B) = a \tag{1}$$

b) In the case that the property of derivability of $\mu$ takes place in any $x\in A$, we say that $\mu$ is *derivable on* $A$.

c) If $\mu$ is derivable on any set $A\in \boldsymbol{U}_n$, then it is called *derivable*.

5.9 The number $a\in \boldsymbol{B}_2$ depending on $x\in A$ and the function $A \ni x \mapsto a\in \boldsymbol{B}_2$ whose existence is stated in 5.8 are called the *derivative of* $\mu$ *in* $x$, respectively the *derivative function of* $\mu$ *in* $x$.

5.10 The derivative of $\mu$ in $x$ and the derivative function of $\mu$ in $x$ are noted with $d\mu(x)$. Other notations are:

- $\frac{d\mu}{dl}(x)$, if $n = 1$
- $\frac{d\mu}{dS}(x)$, if $n = 2$
- $\frac{d\mu}{dV}(x)$, if $n = 3$

5.11 **Remark** The set $B\in \boldsymbol{U}_n$ formed by one element, $x\in A$

$$B = \{x\} \tag{1}$$

has the property that for any $\varepsilon > 0$,

$$x\in B \text{ and } d(B) = 0 < \varepsilon \tag{2}$$

from where it is inferred that, if $\mu$ is derivable in $x$, then $d\mu(x)$, that generally does not depend on $B$, is given by:

$$d\mu(x) = \mu(\{x\}) \qquad (3)$$

5.12 a) We suppose that $\mu$ is a derivable measure on $A$. The set

$$supp_A\ d\mu = \{x \mid x \in A, d\mu(x) = 1\} = \{x \mid x \in A, \mu(\{x\}) = 1\} \qquad (1)$$

is called the *support of* $d\mu$ *on* $A$.

b) If $\mu$ is derivable (on any set $A \in \boldsymbol{U}_n$), then by definition the set

$$supp\ d\mu = \{x \mid x \in \boldsymbol{R}^n, d\mu(x) = 1\} = \{x \mid x \in \boldsymbol{R}^n, \mu(\{x\}) = 1\} \qquad (2)$$

is called the *support of* $d\mu$ (on $\boldsymbol{R}^n$).

5.13 **Theorem** We consider the derivable measure $\mu : \boldsymbol{U}_n \to \boldsymbol{B}_2$ on the closed set $A \in \boldsymbol{U}_n$ ($A$ is compact). Then the set $supp_A\ d\mu$ is finite.

**Proof** Let us suppose that $supp_A\ d\mu$ is infinite, in contradiction with the conclusion of the theorem. Because $A$ is bounded, there exists (Cesaro) a convergent sequence $x^p \in supp_A\ d\mu,\ p \in \boldsymbol{N}$ and the fact that $A$ is closed implies that

$$x = \lim_{p \to \infty} x^p \qquad (1)$$

belongs to $A$. We apply the hypothesis of derivability of $\mu$ in $x$:

$$\exists \varepsilon > 0,\ \forall B \in \boldsymbol{U}_n, (x \in B \text{ and } d(B) < \varepsilon) \Rightarrow \mu(B) = \mu(\{x\}) \qquad (2)$$

We fix $\varepsilon, B$ like above so that for some $x^p \neq x$ it is true in addition $x^p \in B$. The set $B - \{x^p\}$ satisfies the same hypothesis like $B$, that is:

$$x \in B - \{x^p\} \text{ and } d(B - \{x^p\}) \le d(B) < \varepsilon \qquad (3)$$

and the conclusion must be the same:

$$\mu(B - \{x^p\}) = \mu(\{x\}) \qquad (4)$$

It is inferred that:

$$\mu(\{x\}) = \mu(B) = \mu((B - \{x^p\}) \vee \{x^p\}) = \qquad (5)$$

$$= \mu(B - \{x^p\}) \oplus \mu(\{x^p\}) = \mu(\{x\}) \oplus 1$$

The last equation is a contradiction, having its origin in our supposition that $supp_A\ d\mu$ is infinite.

5.14 **Corollary** Let the measure $\mu : \boldsymbol{U}_n \to \boldsymbol{B}_2$.

a) If $\mu$ is derivable on the topological closure $\overline{A}$ of the set $A \in \boldsymbol{U}_n$, then:

a.1) the set $supp_A\ d\mu$ is finite

a.2) $$\forall x \in A, \exists \varepsilon > 0,\ \forall B \in \boldsymbol{U}_n, (x \in B \text{ and } d(B) < \varepsilon) \Rightarrow \mu(B - \{x\}) = 0 \qquad (1)$$

a.3) $$\mu(A) = \underset{x \in A}{\Xi}\ \mu(\{x\}) \qquad (2)$$

a.4) For any partition $A_i \subset A, i \in I$, we have

$$\mu(A) = \underset{i \in I}{\Xi}\ \mu(A_i) \qquad (3)$$

b) If the measure $\mu$ is derivable, then the set $supp\ d\mu$ is locally finite.

**Proof** a.3) If $supp_A\, d\mu$ is empty, then for any $x \in A$ we have that $x \notin supp_A\, d\mu$ and by replacing in 5.8 (1) $B$ with $A$ and $a$ with $\mu(\{x\})$, it results

$$\mu(A) = \mu(\{x\}) = 0 \tag{4}$$

making the statement of the theorem obvious.

We suppose now that

$$supp_A\, d\mu = \{x^1,...,x^p\},\ p \geq 1 \tag{5}$$

There exists a partition $A_1,...,A_p \in \boldsymbol{U}_n$ of $A$ with the property that $x^i \in A_i, i = \overline{1,p}$ and moreover

$$\mu(A) = \mu(\bigvee_{i=1}^{p} A_i) = \Xi_{i=1}^{p} \mu(A_i) = \Xi_{i=1}^{p} \mu(\{x^i\}) = \Xi_{x \in A} \mu(\{x\})(= \pi(p)) \tag{6}$$

b) $\mu$ is derivable on the compacts $\overline{A} \in \boldsymbol{U}_n$ and from a.1) all the sets

$$A \wedge supp\ d\mu = supp_A\ d\mu \tag{7}$$

are finite.

5.15 Let us suppose that $\mu : \boldsymbol{U}_n \to \boldsymbol{B}_2$ is derivable on the topological closure $\overline{A}$ of $A \in \boldsymbol{U}_n$. The binary number

$$\mu(A) = \Xi_{x \in A} \mu(\{x\}) = \Xi_{x \in A} d\mu(x) = \Xi_{x \in \boldsymbol{R}^n} f(x) \cdot d\mu(x) \tag{1}$$

is noted with $\int_A d\mu,\ \int_{\boldsymbol{R}^n} f \cdot d\mu$ or $\int f \cdot d\mu$ and is called the *integral of* $f : \boldsymbol{R}^n \to \boldsymbol{B}_2$ *relative to* $\mu$, where the relation between $f$ and $A$ is by definition the following:

$$A = \{x \mid x \in \boldsymbol{R}^n, f(x) = 1\} = supp\ f \tag{2}$$

5.16 We note with

$$I_{Loc}^{(n)} = \{f \mid f : \boldsymbol{R}^n \to \boldsymbol{B}_2, supp\ f \in Loc_f^{(n)}\} \tag{1}$$

the $\boldsymbol{B}_2$-algebra of the functions with locally finite support, that are called *locally integrable functions.*

5.17 **Theorem** a) The function $g \in I_{Loc}^{(n)}$ defines a derivable measure $\mu^g : \boldsymbol{U}_n \to \boldsymbol{B}_2$ by the formula:

$$\mu^g(A) = \pi(|A \wedge supp\ g|), A \in \boldsymbol{U}_n \tag{1}$$

It is true the relation

$$d\mu^g(x) = g(x), x \in \boldsymbol{R}^n \tag{2}$$

b) Conversely, if $\mu : \boldsymbol{U}_n \to \boldsymbol{B}_2$ is a derivable measure, then there exists in a unique manner the function $g \in I_{Loc}^{(n)}$ so that

$$\mu(A) = \pi(|A \wedge supp\ g|), A \in \boldsymbol{U}_n \tag{3}$$

being also true the relation

$$d\mu(x) = g(x), x \in \boldsymbol{R}^n \tag{4}$$

**Proof** a) The fact that $\mu^g$ is a measure was already proved at 5.6, if we put

$$\mu^{g}(A)=\mu_{supp\ g}(A),\ A\in \boldsymbol{U}_n \qquad (5)$$

and (2) results from

$$d\mu^{g}(x)=\mu^{g}(\{x\})=\pi(|\{x\}\wedge supp\ g\,|)=\begin{cases}\pi(1)=1, if\ \ x\in supp\ g\\ \pi(0)=0, if\ \ x\notin supp\ g\end{cases}=g(x) \qquad (6)$$

b) If $\mu$ is derivable, then $supp\ d\mu\in Loc_{f}^{(n)}$ from 5.14 b) and the function $g:\boldsymbol{R}^n\to \boldsymbol{B}_2$ defined by:

$$g(x)=d\mu(x)=\mu(\{x\}),\ x\in \boldsymbol{R}^n \qquad (7)$$

is locally integrable. As (4) was proved at (7), we prove (3) by taking into account 5.14 a.3):

$$\mu(A)=\underset{x\in A}{\Xi}\ \mu(\{x\})=\pi(|\,A\wedge supp\ d\mu\,|)=\pi(|\,A\wedge supp\ g\,|),\ A\in \boldsymbol{U}_n \qquad (8)$$

5.18 **Corollary** For $g\in I_{Loc}^{(n)}$ and $A\in \boldsymbol{U}_n$ it is defined the integral

$$\int_A g=\int_A d\mu^{g}=\mu^{g}(A)=\pi(|\,A\wedge supp\ g\,|) \qquad (1)$$

## 6. The Lebesgue-Stieltjes Measure

6.1 We say about the function $f:\boldsymbol{R}\to \boldsymbol{B}_2$ that

a) it has a *left limit* in any point $t\in \boldsymbol{R}\vee\{\infty\}$, if (see 3.10)

$$\forall t\in \boldsymbol{R}\vee\{\infty\},\exists t'<t,\exists f(t-0)\in \boldsymbol{B}_2,\forall \xi\in(t',t),\ f(\xi)=f(t-0) \qquad (1)$$

b) it is *left continuous* in any point $t\in \boldsymbol{R}$, if a) is true in any $t\in \boldsymbol{R}$ and moreover:

$$\forall t\in \boldsymbol{R},\ f(t)=f(t-0) \qquad (2)$$

6.2 We fix a function $f$ satisfying the properties from 6.1. We prolong $f$ to $\boldsymbol{R}\vee\{\infty\}$ by left continuity in the point $\infty$:

$$f(\infty)=f(\infty-0) \qquad (1)$$

and we note this new function with $f$ too.

6.3 The relation

$$\mu([[a_1,b_1))\,\Delta\,...\,\Delta\,[[a_n,b_n)))=f(a_1)\oplus f(b_1)\oplus\ ...\ \oplus f(a_n)\oplus f(b_n) \qquad (1)$$

where $a_1,...,a_n,b_1,...,b_n\in \boldsymbol{R}\vee\{\infty\}$ obviously defines an additive function $\mu:Sym^{-}\to \boldsymbol{B}_2$ (see 3.11 for the definition of the symmetrical intervals and of $Sym^{-}$). Our purpose is that of proving the next:

6.4 **Theorem** $\mu$ is a measure.

**Proof** Let $A_n\in Sym^{-}, n\in \boldsymbol{N}$ a sequence of sets that are disjoint two by two with the property that the reunion

$$A=\underset{n\in \boldsymbol{N}}{\vee}A_n \qquad (1)$$

belongs to $Sym^{-}$. We can suppose without loss that all the sets $A_n$ are non-empty.

Case a) $(t_n)_n$ is a real strictly increasing sequence that converges to $b$ and we have:

$$-\infty < a = t_0 < t_1 < t_2 < ... < b \le \infty \quad (2)$$

$$A_n = [t_n, t_{n+1}), n \in \mathbf{N} \quad (3)$$

$$A = [a, b) \quad (4)$$

There exists an $N \in \mathbf{N}$ with

$$n > N \Rightarrow \mu(A_n) = f(t_n) \oplus f(t_{n+1}) = f(b-0) \oplus f(b-0) = 0 \quad (5)$$

and we can write that

$$\underset{n\in \mathbf{N}}{\Xi} \mu(A_n) = \overset{N}{\underset{n=0}{\Xi}} \mu(A_n) = f(t_0) \oplus f(t_1) \oplus f(t_1) \oplus f(t_2) \oplus ... \oplus f(t_N) \oplus f(t_{N+1}) \quad (6)$$

$$= f(t_0) \oplus f(t_{N+1}) = f(a) \oplus f(b-0) = f(a) \oplus f(b) = \mu(A)$$

Case b) $A$ is of the general form

$$A = [a_1, b_1) \vee [a_2, b_2) \vee ... \vee [a_k, b_k) \quad (7)$$

where

$$-\infty < a_1 < b_1 \le a_2 < b_2 \le ... \le a_k < b_k \le \infty \quad (8)$$

We note

$$A_{n,i} = A_n \wedge [a_i, b_i) \quad (9)$$

where $n \in \mathbf{N}$ and $i = \overline{1,k}$ ; we have:

$$\underset{n\in \mathbf{N}}{\Xi} \mu(A_n) = \underset{n\in \mathbf{N}}{\Xi} (\mu(A_{n,1} \vee A_{n,2} \vee ... \vee A_{n,k})) = \quad (10)$$

$$= \underset{n\in \mathbf{N}}{\Xi} (\mu(A_{n,1} \oplus \mu(A_{n,2}) \oplus ... \oplus \mu(A_{n,k})) =$$

$$= \underset{n\in \mathbf{N}}{\Xi} \mu(A_{n,1} \oplus \underset{n\in \mathbf{N}}{\Xi} \mu(A_{n,2}) \oplus ... \oplus \underset{n\in \mathbf{N}}{\Xi} \mu(A_{n,k}) =$$

$$= \mu(\underset{n\in \mathbf{N}}{\vee} A_{n,1}) \oplus \mu(\underset{n\in \mathbf{N}}{\vee} A_{n,2}) \oplus ... \oplus \mu(\underset{n\in \mathbf{N}}{\vee} A_{n,k}) = \quad \text{(from a))}$$

$$= \mu([a_1, b_1)) \oplus \mu([a_2, b_2)) \oplus ... \oplus \mu([a_k, b_k)) =$$

$$= \mu([a_1, b_1) \vee [a_2, b_2) \vee ... \vee [a_k, b_k)) = \mu(A)$$

6.5 The measure $\mu$ that was defined at 6.3 is called the (*left*) *Lebesgue-Stieltjes measure associated to* $f$ .

6.6 The right dual construction is made starting from an $\mathbf{R} \to \mathbf{B}_2$ function (see 6.1) with a right limit in any $t \in \{-\infty\} \vee \mathbf{R}$, right continuous in any $t \in \mathbf{R}$ that is prolonged (see 6.2) to $\{-\infty\} \vee \mathbf{R}$ by right continuity in the point $-\infty$. It is defined then (see 6.3) a measure $Sym^+ \to \mathbf{B}_2$, where $Sym^+$, the dual of $Sym^-$, is the set ring generated by the symmetrical intervals

$$((a,b]] = \begin{cases} (a,b], a < b \\ (b,a], b < a \\ \emptyset \ , a = b \end{cases} \quad (1)$$

where $a, b \in \{-\infty\} \vee \mathbf{R}$ .

6.7 **Theorem** Let $\mu_1 : Sym^- \to \mathbf{B}_2$ an arbitrary measure.

a) The function

$$g(t) = \mu_1([[a,t))) \quad (1)$$

where $a, t \in \mathbf{R} \vee \{\infty\}$ is left continuous on $\mathbf{R} \vee \{\infty\}$.

b) $\mu_1$ is the left Lebesgue-Stieltjes measure associated to $g$.

**Proof** a) It is considered the sequence $(t_n)_n$

$$-\infty < a = t_0 < t_1 < t_2 < \ldots < t \leq \infty \tag{2}$$

that is strictly increasing and convergent to $t$. The sets

$$A_n = [t_n, t_{n+1}), n \in \mathbf{N} \tag{3}$$

belong to $Sym^-$ and are disjoint two by two and their reunion

$$\underset{n \in \mathbf{N}}{\vee} A_n = [a, t) \tag{4}$$

is an element from $Sym^-$ too. It results that there exists $N \in \mathbf{N}$ with

$$n > N \Rightarrow \mu_1(A_n) = \mu_1([t_n, t_{n+1})) = \mu_1([a, t_n)) \oplus \mu_1([a, t_{n+1})) = g(t_n) \oplus g(t_{n+1}) = 0 \tag{5}$$

showing the existence of $g(t-0)$. But

$$g(t) = \mu_1([a,t)) = \mu_1(\underset{n \in \mathbf{N}}{\vee} [t_n, t_{n+1})) = \underset{n \in \mathbf{N}}{\Xi} \mu_1([t_n, t_{n+1})) = \tag{6}$$

$$= \underset{n=0}{\overset{N}{\Xi}} \mu_1([t_n, t_{n+1})) = \underset{n=0}{\overset{N}{\Xi}} (g(t_n) \oplus g(t_{n+1})) = g(t_0) \oplus g(t_{N+1})$$

From the fact that

$$g(t_0) = \mu_1([a,a)) = \mu_1(\emptyset) = 0 \tag{7}$$

$$g(t_{N+1}) = g(t-0) \tag{8}$$

it results, as $t$ is arbitrary, the statement of the theorem.

b) We have that

$$\mu_1([[a_1, b_1)) \Delta \ldots \Delta [[a_n, b_n))) = \tag{9}$$

$$= \mu_1([[a, a_1)) \Delta [[a, b_1)) \Delta \ldots \Delta [[a, a_n)) \Delta [[a, b_n))) =$$

$$= \mu_1([[a, a_1))) \oplus \mu_1([[a, b_1))) \oplus \ldots \oplus \mu_1([[a, a_n))) \oplus \mu_1([[a, b_n))) =$$

$$= g(a_1) \oplus g(b_1) \oplus \ldots \oplus g(a_n) \oplus g(b_n)$$

is true for any $a_1, \ldots, a_n, b_1, \ldots, b_n \in \mathbf{R} \vee \{\infty\}$.

## 7. Measurable Spaces and Measurable Functions. The Integration of the Binary Functions Relative to a Measure

7.1 It is called *measurable space* a pair $(X, \mathbf{U})$ where $X$ is a set and $\mathbf{U} \subset 2^X$ is a ring of subsets of $X$. The sets $A \in \mathbf{U}$ are called *measurable.*

7.2 We say that we have defined a *measurable*, or *integrable function* $f : (X, \mathbf{U}) \to \mathbf{B}_2$ where $(X, \mathbf{U})$ is a measurable space, if it is given the function $f : X \to \mathbf{B}_2$ with $supp\ f \in \mathbf{U}$, i.e. the support of $f$ is a measurable set.

7.3 Recall that the set of the binary functions

$$\mathbf{U}' = \{f \mid f : (X, \mathbf{U}) \to \mathbf{B}_2,\ f \textit{ is measurable}\} \tag{1}$$

(see 1.13, 1.14) is a $\mathbf{B}_2$-algebra relative to the obvious laws. As ring, it is isomorphic with $\mathbf{U}$.

7.4 Let $(X, \mathbf{U})$ a measurable space and $M \subset X$. Because the set

$$\boldsymbol{U}\wedge M=\{A\wedge M\mid A\in\boldsymbol{U}\}\qquad(1)$$
is a set ring, the pair $(M,\boldsymbol{U}\wedge M)$ is a measurable space, called *measurable subspace* of $(X,\boldsymbol{U})$.

7.5 **Proposition** a) If $f:(X,\boldsymbol{U})\to\boldsymbol{B}_2$ is measurable, then its restriction $f_{|M}:(M,\boldsymbol{U}\wedge M)\to\boldsymbol{B}_2$ is measurable.

b) If $g:(M,\boldsymbol{U}\wedge M)\to\boldsymbol{B}_2$ is measurable, then it can be prolonged to a measurable function $f:(X,\boldsymbol{U})\to\boldsymbol{B}_2$.

7.6 Let us suppose that $(X,\boldsymbol{U})$ is a measurable space, $f:(X,\boldsymbol{U})\to\boldsymbol{B}_2$ is a measurable function and $\mu:\boldsymbol{U}\to\boldsymbol{B}_2$ is a measure. The number $\mu(supp\ f)$ is called the *integral of* $f$ *relative to* $\mu$ and is noted with $\int f\cdot d\mu$.

7.7 Let $f_n, f\in\boldsymbol{U}', n\in\boldsymbol{N}$.

a) If
$$supp\ f_0\subset supp\ f_1\subset supp\ f_2\subset...\qquad(1)$$
$$\bigvee_{n\in\boldsymbol{N}} supp\ f_n=supp\ f\qquad(2)$$
then we say that $f_n$ *converges*, or *tends increasingly to* $f$ and this fact is noted with $f_n\uparrow f$.

b) If
$$supp\ f_0\supset supp\ f_1\supset supp\ f_2\supset...\qquad(1)$$
$$\bigwedge_{n\in\boldsymbol{N}} supp\ f_n=supp\ f\qquad(2)$$
then we say that $f_n$ *converges*, or *tends decreasingly to* $f$ and this fact is noted with $f_n\downarrow f$.

c) In one of the situations from a), b) we say that $f_n$ *converges*, or *tends monotonously to* $f$ and the notation is $f_n\updownarrow f$.

7.8 Let us suppose that $f,g:(X,\boldsymbol{U})\to\boldsymbol{B}_2$ are measurable and that $\mu:\boldsymbol{U}\to\boldsymbol{B}_2$ is a measure. We say that $f$ and $g$ are *equal almost everywhere* and we write this fact with
$$f=g\ \text{ a.e.}\qquad(1)$$
if
$$\mu(\{x\mid f(x)\neq g(x)\})=\mu(supp\ f\ \Delta\ supp\ g)=0\qquad(2)$$
or, in an equivalent manner, if
$$\mu(supp\ f)=\mu(supp\ g)\qquad(3)$$

7.9 **Proposition** The function $\boldsymbol{U}'\ni f\mapsto\int f\cdot d\mu\in\boldsymbol{B}_2$ satisfies the following properties:

a) it is linear

b) $f_n\updownarrow f\Rightarrow\int f_n\cdot d\mu\to\int f\cdot d\mu$

c) $f=g$ a.e. $\Leftrightarrow\int f\cdot d\mu=\int g\cdot d\mu$

where $f_n,f,g\in\boldsymbol{U}', n\in\boldsymbol{N}$.

**Proof** b) is a restatement of 4.2 a) and 4.3 a).

7.10 **Corollary** If $f_n \in U', n \in N$ converges to $0$ decreasingly, then

$$\int f_n \cdot d\mu \to 0 \qquad (1)$$

7.11 Let $f_n, f : (X,U) \to B_2, n \in N$ measurable and $\mu : U \to B_2$ a measure. We say that $f_n$ *tends to* $f$ *in measure* and we note this property with $f_n \underset{\mu}{\to} f$ if

$$\int f_n \cdot d\mu \to \int f \cdot d\mu \qquad (2)$$

7.12 Let $f, \chi_A : X \to B_2$ two functions, where $f$ is arbitrary and $\chi_A$ is the characteristic function of the set $A \subset X$ . If $A \wedge supp\ f \in U$ - condition that is called *of integrability*- then the number

$$\int_A f \cdot d\mu \overset{def}{=} \int (\chi_A \cdot f) \cdot d\mu \qquad (1)$$

is called the *integral of* $f$ *, on* $A$*, relative to* $\mu$.

7.13 The function $f \cdot \mu : U \to B_2$ defined by

$$(f \cdot \mu)(A) = \int_A f \cdot d\mu \qquad (1)$$

where $A$, $supp\ f \in U$ is a measure, that coincides with the restriction of $\mu$ at $supp\ f$ .

## 8. Riemann Integrals

8.1 We end the paper with a short paragraph that introduces the Riemann integrals of the $f : R \to B_2$ functions (generalizations are possible to $f : R^n \to B_2$ functions). The main feature for this type of integral is considering the set ring $R_f(R)$ and the finite Boolean measure (see 3.8) $\mu_f^R : R_f(R) \to B_2$.

8.2 For the set $A \subset R$ , the property $A \wedge supp\ f \in R_f(R)$ (see 7.12) is called the condition of *Riemann integrability of* $f$ *on* $A$. If it is fulfilled, we say that $f$ is *Riemann integrable,* or *integrable in the sense of Riemann, on* $A$.

8.3 **Special cases** for 8.2. a) $[[a,b)) \wedge supp\ f \in R_f(R)$ (see 3.12 (1) for the definition of $I_{[[a,b))}$), $a,b \in R \vee \{\infty\}$. These functions are called *left integrable* (in the sense of Riemann) *from* $a$ *to* $b$.

b) $R \wedge supp\ f \in R_f(R)$ (see 3.12 (2) for the definition of $I_\infty$). These functions are called *integrable* (in the sense of Riemann).

c) $\forall a,b \in R, (a,b) \wedge supp\ f \in R_f(R)$ (see 5.4, 5.5, 5.16 for the definition of $I_{Loc}^{(1)}$). These functions are called *locally integrable* (in the sense of Riemann) and they have a locally finite support.

d) $\forall a,b \in R \vee \{\infty\}, (a,b) \wedge supp\ f \in R_f(R)$ defines the $B_2$-algebra of functions $I_{Sup}$. We say about these functions that they are *left integrable* (in the sense of Riemann) and that they have the support *superiorly finite*, dual notion to that of inferiorly finite set that was defined at 3.9.

8.4 If $f$ is Riemann integrable on $A$, then the number (see 7.12 (1))

$$\int_A f \cdot d\mu_f^{\boldsymbol{R}} = \mu_f^{\boldsymbol{R}}(A \wedge supp\ f) = \mathop{\Xi}_{x \in A} f(x) \qquad (1)$$

is called the *integral, in the sense of Riemann, of* $f$ *, on* $A$.

8.5 **Special cases** for 8.4 a) $f \in I_{[[a,b))}$, $a,b \in \boldsymbol{R} \vee \{\infty\}$; the integral $\int_{[[a,b))} f \cdot d\mu_f^{\boldsymbol{R}}$ is noted with $\int_a^{b} \bar{} \ f$ and is called the *left integral* (in the sense of Riemann) *of* $f$ *from* $a$ *to* $b$.

b) $f \in I_\infty$; the integral $\int_{\boldsymbol{R}} f \cdot d\mu_f^{\boldsymbol{R}}$ is usually noted with $\int_{-\infty}^{\infty} f$ and is called the *integral* (in the sense of Riemann) *of* $f$.

8.6 The cases 8.3 a) and 8.5 a) have right duals, that refer to symmetrical intervals of the form $((a,b]], a,b \in \{-\infty\} \vee \boldsymbol{R}$ (see 6.6).

8.7 We define the subring of sets $Sym' \subset Sym^-$ to be the one that is generated by the symmetrical intervals $[[a,b)), a,b \in \boldsymbol{R}$ (at $Sym^-$ we had $[[a,b)), a,b \in \boldsymbol{R} \vee \{\infty\}$).

8.8 a) Let us suppose that $f \in I_{Loc}^{(1)}$. Then the measure $f \cdot \mu_f^{\boldsymbol{R}} : Sym' \to \boldsymbol{B}_2$ (see 7.13) is called the *indefinite integral* of $f$.

b) The function $F^- : \boldsymbol{R} \to \boldsymbol{B}_2$, which is defined in the next manner:

$$F^-(t) = f \cdot \mu_f^{\boldsymbol{R}}([[a,t))), t \in \boldsymbol{R} \qquad (1)$$

where $a \in \boldsymbol{R}$ is a parameter is called the *left primitive* of $f$.

c) The left primitive $F^-(t)$ has a left limit and it is left continuous in any $t \in \boldsymbol{R}$.

8.9 If at 8.8 $f \in I_{Sup}$ (where $I_{Sup} \subset I_{Loc}^{(1)}$), then $f \cdot \mu_f^{\boldsymbol{R}}$ is extended to $Sym^-$ and $F^-$ is extended to $\boldsymbol{R} \vee \{\infty\}$, by left continuity in the point $\infty$. $f \cdot \mu_f^{\boldsymbol{R}}$ is in this situation the left Lebesgue-Stieltjes measure associated to $F^-$ (see 6.3).

8.10 Together with the duals in the left-right sense that have appeared having their origin in the order of $\boldsymbol{R}$, the previous notions have also another type of duality, so called in the algebraical sense, resulting by the replacement of $0$ with $1$ and viceversa, to be compared, from the table 1.1, the laws '$\oplus$' and '$\otimes$'. For example, the algebraical dual of $\int_a^{b} \bar{} \ f$ is defined like this:

$$\int_a^{b} {}^{*-} f = \bigotimes_{x \in [[a,b))} f(x) \qquad (1)$$

**Bibliography**

[1] Nicu Boboc, Gheorghe Bucur, Masura si capacitate (measure and Capacity, in Romanian), Ed. stiintifica si enciclopedica, Bucuresti, 1985

[2] Miron Nicolescu, Analiza matematica (Mathematical Analysis, in Romanian), vol 1-3, Ed. Tehnica, Bucuresti, 1957